\documentclass[12pt]{amsart}
\font\emailfont=cmtt10

\usepackage{amsmath,amsthm,amsfonts,amscd,flafter,epsf}
\usepackage[small,nohug,heads=littlevee]{diagrams} 
\diagramstyle[labelstyle=\scriptstyle] 

\newcommand{\HF}{HF}

\newtheorem{theorem}{Theorem}[section]

\newtheorem{lemma}[theorem]{Lemma}

\newtheorem{defn}[theorem]{Definition}

\def\endproof{\relax\ifmmode\expandafter\endproofmath\else
  \unskip\nobreak\hfil\penalty50\hskip.75em\hbox{}\nobreak\hfil\bull
  {\parfillskip=0pt \finalhyphendemerits=0 \bigbreak}\fi}
\def\endproofmath$${\eqno\bull$$\bigbreak}
\def\bull{\vbox{\hrule\hbox{\vrule\kern3pt\vbox{\kern6pt}\kern3pt\vrule}\hrule}}

\newcommand{\R}{\mathbb{R}}

\newcommand{\Z}{\mathbb{Z}}

\newcommand{\cm}{\cdot}

\newcommand{\ModSWfour}{\mathcal{M}}
\newcommand{\ModFlow}{\ModSWfour}

\newcommand\abuts\Rightarrow
\newcommand\Sym{\mathrm{Sym}}

\newcommand\relspinc{\underline{\spinc}}

\newcommand\x{\mathbf x}

\newcommand\q{\mathbf q}

\newcommand\y{\mathbf y}

\newcommand\ModSphere{\ModFlow\left({\mathbb S}\longrightarrow 
\Sym^{g-1}(\Sigma_{1})\times \Sym^2(\Sigma_{2})\right)}
\newcommand\ModSpheres\ModSphere

\newcommand\HFm{\HF^-}

\newcommand\Mas{\mu}
\newcommand\UnparModSp{\widehat \ModSp}
\newcommand\UnparModFlow\UnparModSp
\newcommand\Mod\ModSp

\newcommand{\spinc}{\mathfrak s}

\newcommand\ModMaps{\mathcal M}
\newcommand\ModSp\ModMaps

\newcommand\Ta{{\mathbb T}_{\alpha}}
\newcommand\Tb{{\mathbb T}_{\beta}}

\newcommand\Torus{\mathbb{T}}

\newcommand\alphas{\mbox{\boldmath$\alpha$}}

\newcommand\betas{\mbox{\boldmath$\beta$}}

\newcommand\Dual{\mathcal D}
\newcommand\Duality\Dual

\newcommand\eD{\widetilde{D}}

\newcommand\Pent{\mathrm{Pent}}
\newcommand\EmptyPent{\Pent^\circ}
\newcommand{\fa}{\widehat f}
\newcommand{\ga}{\widehat g}
\newcommand{\ha}{\widehat h}

\newcommand{\fm}{f^-}
\newcommand{\gm}{g^-}
\newcommand{\hm}{h^-}

\newcommand\Kp{\Knot_+}
\newcommand\Km{\Knot_-}

\newcommand\oLink{\vec{L}}
\newcommand\oKnot{\vec{K}}

\newcommand\Gen{\mathfrak{S}}
\def\Interior{\mathrm{Int}}
\newcommand\EmptyRect{\Rect^\circ}
\newcommand\Cm{C^-}

\newcommand\dm{\partial^-}

\newcommand\Rect{\mathrm{Rect}}

\newcommand\Os{\mathbb O}

\newcommand\Xs{\mathbb X}

\newcommand\Knot{\mathcal K}
\newcommand\Link{\mathcal L}

\newcommand\CFLm{\mathrm{CFL}^-}
\newcommand\HFLm{\mathrm{HFL}^-}

\newcommand\spincrel\relspinc

\newcommand\CFK{\mathrm{CFK}}
\newcommand\HFK{\mathrm{HFK}}

\newcommand\CFKa{\widehat\CFK}
\newcommand\CFKm{\CFK^-}

\newcommand\HFKa{\widehat\HFK}
\newcommand\HFKm{\HFK^-}

\title[{On the skein exact squence for knot Floer homology}]
{On the skein exact squence for knot Floer homology}

\author[Peter Ozsv{\'a}th]{Peter Ozsv\'ath}
\address{Department of
Mathematics, Columbia University, 
New York, NY 10027 \newline
\indent{\emailfont{petero@math.columbia.edu}}}
\thanks{PSO was supported by NSF grant number DMS-0505811 and  FRG-0244663}

\author[Zolt{\'a}n Szab{\'o}]{Zolt{\'a}n Szab{\'o}} 
\address{Department of
Mathematics, Princeton University, New Jersey 08544 \newline
\indent{\emailfont{szabo@math.princeton.edu}}}
\thanks{ZSz was supported by NSF grant number DMS-0406155 and  FRG-0244663}

\begin{document}
\begin{abstract}
  The aim of this paper is to study the skein exact sequence for knot
  Floer homology. We prove precise graded version of this sequence,
  and also one using $\HFm$. Moreover, a complete argument is also given
  purely within the realm of grid diagrams.
\end{abstract}

\maketitle

\section{Introduction}

Knot Floer homology is an invariant for knots in $S^3$ defined using
Heegaard diagrams and holomorphic disks~\cite{Knots},
\cite{RasmussenThesis}. This invariant can be used to construct a
bigraded group $\HFKa$, endowed with an {\em Alexander} and a {\em
  Maslov} grading, has as its Euler characteristic the Alexander
polynomial of the knot. Another variant gives a bigraded Abelian group
$\HFm$, which is a module over the polynomial ring $\Z[U]$, and whose
specialization (in a suitable sense) to $U=0$ gives $\HFKa$.

\begin{figure}[ht]
\mbox{\vbox{\epsfbox{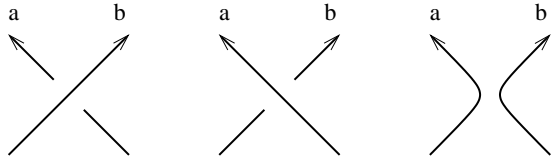}}}
\caption{\label{fig:Crossings}
{\bf{Skein triple.}} Diagram of a positive crossing, a negative crossing, and the
(oriented) resolution respectively.}
\end{figure}

The traditional skein relation for the Alexander polynomial translates
into this context into a long exact sequence which relates $\HFKa$ of
a knot with a distinguished positive crossing $\Knot_+$, the knot
Floer homology of the oriented resolution $\Knot_0$ of that crossing
(which is a link), and also the knot Floer homology of the knot
$\Knot_-$ obtained by changing the distinguished positive crossing in
$\Knot_+$ to a negative crossing, see Figure~\ref{fig:Crossings}. The first version of this exact
triangle appeared in~\cite{Knots}, where the term involving $\Knot_0$
is defined using a suitable generalization of knot Floer homology to
an invariant for oriented links, which we denote by $\HFKa(\Knot_0)$.

Link Floer homology is given a more general definition
in~\cite{Links}, as a multi-graded theory whose Euler characteristic
is the multi-variable Alexander polynomial. Algebraically, the
invariant $\HFLm(\Link)$ is a multi-graded theory which is the
homology of a chain complex $\CFLm(\Link)$ over $\Z[U_1,...,U_\ell]$,
where the formal variables $U_i$ are in one-to-one correspondence with
the components of the link. The invariant appearing in the earlier
skein exact sequence is the homology group $\HFKa(\Link)$ gotten by
setting all the $U_i=0$, and adding up all of the ``Alexander
gradings''.

In~\cite{MOS}, link Floer homology is given a purely combinatorial
calculation via ``grid diagrams''. This thread is pursued further
in~\cite{MOST}, where the basics of the theory are developed from a
purely combinatorial point of view.

The aim of the present paper is to give a different proof of the skein
exact sequence for knot Floer homology. The advantages of this proof
is that it generalizes to the case of $\HFKm$, and also we can give
more precise grading information about the maps. Moreover, this
perspective can be applied readily to give another (quite similar)
proof which works purely within the context of grid diagrams. Aside
from an aesthetic benefit, this also gives a direct combinatorial way
to calculate the maps appearing in the skein exact sequence.

\begin{theorem}
  \label{thm:GradedSkein}
  Let $\Knot_+$, $\Knot_0$, and $\Knot_-$ be three links, which differ
  at a single crossing as indicated by the notation.  Then, if the two
  strands meeting at the distinguished crossing in $\Knot_+$ belong to
  the same component, so that in the oriented resolution the two
  strands corresponding to two distinct components $a$ and $b$ of
  $\Knot_0$, then there are long exact sequences
  \begin{tiny}
  \[
  \begin{CD}
    ...@>>>\HFKa_m(\Knot_+,s)@>{\fa}>>\HFKa_{m}(\Knot_-,s)@>{\ga}>>\HFKa_{m-1}(\Knot_0,s)@>{\ha}>>\HFKa_{m-1}(\Knot_+,s)@>{\fa}>>... \\
    ...@>>>\HFKm_m(\Knot_+,s)@>{\fm}>>\HFKm_m(\Knot_-,s)@>{\gm}>>H_{m-1}\left(\frac{\CFLm(\Knot_0)}{U_a-U_b},s\right)@>{\hm}>>\HFK_{m-1}(\Knot_+,s)@>{\fa}>>...
  \end{CD}
  \]
  \end{tiny}
  If they belong to different components, we have a long exact sequence
  \begin{tiny}
  \[
  \begin{CD}
    ...@>>>\HFKa_{m}(\Knot_+,s)@>{\fa}>>\HFKa_{m}(\Knot_-,s)@>{\ga}>>\left(\HFKa(\Knot_0)\otimes V\right)_{m-1,s} @>{\ha}>>\HFKa_{m-1}(\Knot_+,s)@>{\fa}>>...\\
    ...@>>>\HFKm_{m}(\Knot_+,s)@>{\fm}>>\HFKm_{m}(\Knot_-,s)@>{\gm}>>\left(\HFKm(\Knot_0)\otimes W\right)_{m-1,s} @>{\hm}>>\HFKa_{m-1}(\Knot_+,s)@>{\fa}>>...,
  \end{CD}
  \]
  \end{tiny}
  where here $V$ is the bigraded module
  $$V_{m,s}\cong \left\{\begin{array}{ll}
      \Z^2 &  {(m,s)=(-1,0)} \\
      \Z & {(m,s)=(0,-1)~\text{or}~(1,0)} \\
      0 &{\text{otherwise,}}\end{array}\right.$$
  and $W$ is the bigraded module 
  $$W_{m,s}\cong \left\{\begin{array}{ll}
      \Z & {(m,s)=(0,0)~\text{or}~(1,1)} \\
      0 &{\text{otherwise.}}\end{array}\right.$$
\end{theorem}

The reader is warned: there are two natural conventions on Maslov
grading, one which takes half-integral values (cf.~\cite{Knots}), and
the other which always takes integral values (cf.~\cite{Links}). In
the above statement, we have adopted the latter convention. 

A version of Theorem~\ref{thm:GradedSkein} appears in~\cite{Knots},
except that the map defined there $f$ is not known to preserve Maslov
gradings. This renders that version of the skein exact sequence
somewhat cumbersome to use. It is interesting to note that the
the module $V$ appears in for quite different reasons in the two approaches.

Two slightly different proofs of Theorem~\ref{thm:GradedSkein} are
given. The first uses pseudo-holomorphic disks. The second is a
combinatorial proof involving grid diagrams. This proof is slightly
more awkward, as one cannot use a fixed grid diagram for all three
knots; and of course, it is slightly less awkward in that it is a
purely combinatorial argument, and the maps can be defined by explicit
counts of polygons. Both proofs can be seen as a double iteration of
the skein relating involving singular knots from~\cite{Resolutions},
defined using Floer homology for singular knots from~\cite{SingLink}.
We have however chosen to give a more self-contained proof of
Theorem~\ref{thm:GradedSkein} making no explicit reference to Floer
homology for singular links; but our proof here is very similar in
spirit to the proof of the skein sequence involving singular
links,~\cite[Theorem~\ref{Resolutions:thm:SkeinExactSequenceIntro}]{Resolutions}.

It is possible that the map defined here $\fa$ differs from the one
used in~\cite{Knots}. It also seems different from the one used
in~\cite{Audoux}.  In the next section, we briefly recall knot Floer
homology, and set up our notation. In Section~\ref{sec:Proof}, we
state and prove a theorem which specializes readily to
Theorem~\ref{thm:GradedSkein}. 

\subsection{Acknowledgements}
We wish to thank Benjamin Audoux, {\'E}tienne Gallais, Matt Hedden,
Ciprian Manolescu, and Dylan Thurston for interesting discussions.

\section{Floer homology of knots and links}
\label{sec:Knots}

Knot Floer homology is a bigraded Abelian group associated to a knot
in $S^3$, cf.~\cite{Knots}, \cite{RasmussenThesis}.  We will briefly
sketch this construction, and refer the reader to the above sources
for more details.

Let $\Sigma$ be a surface of genus $g$, let
$\alphas=\{\alpha_1,...,\alpha_{g+n-1}\}$ be a collection of pairwise
disjoint, embedded closed curves in $\Sigma$ which span a
$g$-dimensional subspace of $H_1(\Sigma)$. This specifies a handlebody
$U_\alpha$ with boundary $\Sigma$. Moreover,
$\alpha_1\cup...\cup\alpha_{g+n-1}$ divides $\Sigma$ into $n$
components, which we label
$$\Sigma-\alpha_1-...-\alpha_{g+n-1}={\mathfrak A}_1\coprod...\coprod {\mathfrak A}_n.$$
Fix another such collection of
curves $\betas=\{\beta_1,...,\beta_{g+n-1}\}$, giving another
handlebody $U_\beta$. 
Write
$$\Sigma-\beta_1-...-\beta_{g+n-1}={\mathfrak B}_1\coprod...\coprod
{\mathfrak B}_n.$$
Let $Y$ be the three-manifold specified by the
Heegaard decomposition specified by the handlebodies $U_\alpha$ and
$U_\beta$.  Choose collections of disjoint points
$\Os=\{O_1,...,O_n\}$ and $\Xs=\{X_1,...,X_n\}$, which are distributed
so that each region ${\mathfrak A}_i$ and ${\mathfrak B}_i$ contains
exactly one of the points in $\Os$ and also exactly one of the points
in $\Xs$.  We can use the points $\Os$ and $\Xs$ to construct an
oriented, embedded one-manifold $\oLink$ in $Y$ by the following
procedure. Let $\xi_i$ be an arc connecting the point in $\Xs\cap
{\mathfrak A}_i$ with the point in $\Os\cap {\mathfrak A}_i$, and let
$\xi_i'$ be its pushoff into $U_\alpha$ i.e. the endpoints of $\xi_i'$
coincide with those of $\xi_i$ (and lie on $\Sigma$), whereas its
interior is an arc in the interior of $U_\alpha$. The arc is endowed
with an orientation, as a path from an element of $\Xs$ to an element
of $\Os$.  Similarly, let $\eta_i$ be an arc connecting $\Os\cap
{\mathfrak B}_i$ to $\Xs\cap {\mathfrak B}_i$, and $\eta_i'$ be its
pushoff into $U_\beta$. Putting together the
$\xi_i'$ and $\eta_i'$, we obtain an oriented link $\oLink$ in $Y$.

\begin{defn}
The data $(\Sigma,\alphas,\betas,\Os,\Xs)$ is called a 
{\em pointed Heegaard diagram compatible with the oriented link
$\oLink\subset Y$}.
\end{defn}

An oriented link in a closed three-manifold $Y$ always admits a
compatible pointed Heegaard diagram. In this article, we will restrict
attention to the case where the ambient three-manifold $Y$ is $S^3$.

Consider now the $g+n-1$-fold symmetric product of the surface
$\Sigma$, $\Sym^{g+n-1}(\Sigma)$.  This space is equipped with a pair
of tori
\begin{eqnarray*}
\Ta=\alpha_1\times...\times\alpha_{g+n-1} 
&{\text{and}}&
\Tb=\beta_1\times...\times\beta_{g+n-1}.
\end{eqnarray*}
Knot Floer homology is defined using a suitable variant of Lagrangian
Floer homology for this pair of subsets.

Specifically, let $\Gen$ denote the set of intersection points
$\Ta\cap\Tb\subset\Sym^{g+n-1}(\Sigma)$.
Let $\CFKm(\oLink)$ be the free module over $\Z[U_1,...,U_n]$
generated by elements of $\Gen$, where here the $\{U_i\}_{i=1}^n$ are
indeterminates.

To construct bigradings, consider functions
\begin{eqnarray*}
A\colon \Gen\times\Gen\longrightarrow \Z
&{\text{and}}&
M\colon \Gen\times\Gen\longrightarrow \Z
\end{eqnarray*}
defined as follows.
Given $\x,\y\in\Gen$, let $$A(\x,\y)=\sum_{i=1}^n
(X_i(\phi)-O_i(\phi)),$$ where $\phi\in\pi_2(\x,\y)$ is any Whitney disk
from $\x$ to $\y$, and $X_i(\phi)$ resp. $O_i(\phi)$ is the algebraic
intersection number of $\phi$ with the submanifold $\{X_i\}\times
\Sym^{g+n-2}(\Sigma)$
resp. $\{O_i\}\times
\Sym^{g+n-2}(\Sigma)$. Also, let
$$M(\x,\y)=\Mas(\phi)-2\sum_{i=1}^n O_i(\phi),$$
where $\Mas(\phi)$
denotes the Maslov index of $\phi$; see~\cite{LipshitzCyl} for an
explicit formula in terms of data on the Heegaard diagram.  Both
$A(\x,\y)$ and $M(\x,\y)$ are independent of the choice of $\phi$ in
their definition. There are functions $A\colon \Gen\longrightarrow \Z$
and $M\colon \Gen\longrightarrow \Z$ both of which are uniquely
specified to overall translation by the formulas
\begin{eqnarray}
  \label{eq:RelToAbsolute}
A(\x)-A(\y)=A(\x,\y)
&{\text{and}}&
M(\x)-M(\y)=M(\x,\y).
\end{eqnarray}
The additive indeterminacy in $A$ and $M$ can be removed, as we explain
at the end of the present subsection.

Let $\CFKm(\oLink)$ be the free module over $\Z[U_1,...,U_n]$
generated by $\Gen$. This module inherits a bigrading from
the functions $M$ and $A$ above, with the additional convention that
multiplication by $U_i$ drops the Maslov grading by two, and the
Alexander grading by one.

We define
the differential $$\partial\colon \CFKm(\oLink)\longrightarrow
\CFKm(\oLink)$$
by the formula:
\begin{equation}
\label{eq:DefPartial}
\partial (\x)=\sum_{\y\in\Gen}\,
\sum_{\left\{\phi\in\pi_2(\x,\y)\big| 
\begin{tiny}
\begin{array}{l}
\Mas(\phi)=1 \\
X_i(\phi)=0~~~~~~~~~ \forall i=1,...,n
\end{array}
\end{tiny}
\right\}}\!\!
\#\UnparModFlow(\phi)\cm U_1^{O_1(\phi)}\cdots U_n^{O_n(\phi)}\cm \y.  
\end{equation}
Here, $\UnparModFlow(\phi)$ denotes the  moduli space of
pseudo-holomorphic disks representing the homotopy class $\phi$,
divided out by the action of $\R$. The signed count $\#\UnparModFlow(\phi)$
is associated to an orientation system $\epsilon$, which counts
boundary degenerations with boundary entirely inside $\Ta$ with multiplicity
$+1$ and those with boundary entirely inside $\Tb$ with multiplicity $-1$.
When $\oLink=\oKnot$ is a knot, it it is sometimes convenient to
consider instead the complex $\CFKa(\oKnot)=\CFKm(\oKnot)/(U_1=0)$.
The homology groups $\HFKm(\oKnot)=H_*(\CFKm(\oKnot))$ and
$\HFKa(\oKnot)=H_*(\CFKa(\oKnot))$ are knot invariants~\cite{Knots},
\cite{RasmussenThesis}, see also~\cite{Links}, \cite{MOS} for the 
case of multiple basepoints, and also~\cite{MOST} for a further discussion
of signs. The bigradings on the complex induce bigradings on the homology
\begin{eqnarray*}
\HFKm(\oKnot)=\bigoplus_{m,s}\HFKm_{m}(\oKnot,s) &{\text{and}}&
\HFKa(\oKnot)=\bigoplus_{m,s}\HFKa_{m}(\oKnot,s).
\end{eqnarray*}

For an $\ell$ component link, we consider $\CFKm(\oLink)$ as a module
over $\Z[U_1,...,U_\ell]$, where there is one variable $U_i$
corresponding to each component of $\oLink$.  
In this case, it is
natural to consider
$\CFKa(\oLink)=\CFKm(\oLink)/\{U_1=...=U_\ell=0\}$, and their associated
bigraded homology modules
\begin{eqnarray*}
\HFKm(\oLink)=\bigoplus_{m,s}\HFKm_{m}(\oLink,s) &{\text{and}}&
\HFKa(\oLink)=\bigoplus_{m,s}\HFKa_{m}(\oLink,s).
\end{eqnarray*}
In fact, $\HFKa(\oLink)$ was first defined in~\cite{Knots} using a
slightly different construction, but the equivalence of the two
constructions was established
in~\cite[Theorem~\ref{Links:thm:IdentifyWithLinkHomology}]{Links}.
(In fact, in~\cite{Links}, a more general multi-graded theory is
defined, with one Alexander grading for each component of the
link. The present Alexander grading can be thought of as the sum of
these $\ell$ Alexander gradings. We will have no need for this more
general construction in the present paper.)

We have defined the bigradings only up to additive constants.  This
indeterminacy can be removed with the following conventions.  Dropping
the condition that all the $X_i(\phi)=0$ in the differential for
$\CFKa$, we obtain another chain complex which retains its Maslov
grading, and whose homology is isomorphic to $\Z$
(cf. \cite[Theorem~\ref{Links:thm:InvarianceHFm}]{Links}). 

Similarly, the Alexander grading can be characterized as follows.  If
we consider the complex $C=\CFKm/(U_1=....U_n=1)$. This complex
retains a $\Z$-grading by $N=M-2A$, and its homology is isomorphic to
$H_*(T^{n-1})$ as a relatively graded Abelian group. We fix the
additive constant in the $N$-grading by the requirement that
$$H_m(C) \cong H_{m+2\ell-n-1}(T^{n-1}).$$ This in turn pins down
that additive indeterminacy of $A$.

\subsection{Grid diagrams}
\label{subsec:Grids}

Knot Floer homology has a combinatorial description for Heegaard
diagrams associated to grid presentations according to~\cite{MOS}, cf.
also~\cite{MOST}, \cite{SarkarWang}. 

A grid diagram $G$ is a Heegaard diagram for a knot, where the Heegaard
surface is a torus, and all the $\alphas$ (and the $\betas$) are
parallel, homologically non-trivial circles. We draw the $\alphas$ as
horizontal, and the $\betas$ as vertical. The only non-trivial
contributions in the differential are given by rectangles (and each
such rectangles counts with a sign and also a product of variables 
associated to the the squares marked $\Os$ inside the rectangle).
See~\cite{MOST} for a development of this complex which is logically
independent of holomorphic curve techniques, including a proof of knot
invariance. 

More precisely, if $\{\alpha_1,...,\alpha_n\}$ denote the horizontal
circles and $\{\beta_1,...,\beta_n\}$ are the vertical ones, our
generating set $\Gen$ consists of permutations $\sigma$, which we
think of as $n$-tuples of intersection points $\x$, $x_i=\alpha_i\cap
\beta_{\sigma(i)}$. There are four embedded rectangles in the torus
whose boundary consists of two segments within the $\alphas$ and two
segments with in the $\betas$, and whose four corners are points from
$\x$ and $\y$. Two of these rectangles are oriented so that their
oriented boundary meets the $\alphas$ in a pair of arcs going from
points in $\x$ to points in $\y$. We say that those are the two
rectangles {\em from $\x$ to $\y$}, and we let $\Rect(\x,\y)$ denote
this set of rectangles. If $r\in\Rect(\x,\y)$ has the property that
its interior contains none of the points from $\x$ or $\y$, then we
say $r$ is an {\em empty rectangle}.

In~\cite{MOST}, we verify the existence of a map $\epsilon\colon
\EmptyRect\longrightarrow \{\pm 1\}$ whith the following properties:
\begin{itemize}
\item if $\x,\y,\y',\q\in\Gen$ are generators and $r_1\in
  \EmptyRect(\x,\y)$, $r_2\in\EmptyRect(\y,\q)$ and $r_1'\in
  \EmptyRect(\x,\y')$, $r_2\in\EmptyRect(\y',\q)$, then
  $\epsilon(r_1)\epsilon(r_2)=-\epsilon(r_1')\epsilon(r_2')$
\item if $r_1\in\EmptyRect(\x,\y)$, $r_2\in\EmptyRect(\y,\x)$ are a
  pair of rectangles whose union forms a vertical annulus, then
  $\epsilon(r_1)\cm \epsilon(r_2)=-1$
\item if $r_1\in\EmptyRect(\x,\y)$, $r_2\in\EmptyRect(\y,\x)$ are a
  pair of rectangles whose union forms a horizontal annulus, then
  $\epsilon(r_1)\cm \epsilon(r_2)=1$
\end{itemize}
The chain complex associated to a grid diagram is freely generated by $\Gen$ over $\Z[U_1,...,U_n]$,
with differential given by
$$\dm(\x) = \sum_{\y\in\Gen}\sum_{\{r\in\EmptyRect(\x,\y)\big| r\cap \Xs=\emptyset\}} \epsilon(r)\cm  U_1^{O_1(r)}\cm...\cm U_n^{O_n(r)} \cm \y.$$
As in~\cite{MOS}, this is a special case of the knot Floer homology chain complex considered earlier.

\section{Proofs of the skein sequence}
\label{sec:Proof}

Theorem~\ref{thm:GradedSkein} follows from the following more general
result, Theorem~\ref{thm:StrongForm}, which we state after introducing
a few preliminaries.

Let $\Knot_+$ be a positive crossing, and label its two outgoing edges by
$a$ and $b$, and its two in-coming ones by $c$ and $d$, so that $b$ and $c$
are connected by the crossing, and $a$ and $c$ are connected in the resolution.

Recall that $U_b-U_c$ is an endomorphism of the chain complex
$R=\CFKm(\Knot_0)$, which drops Alexander grading by one and Maslov
grading by two.  Thus, we can form its mapping cone, which is a
bigraded chain complex defined as follows.  Letting $R_{s,d}$ denote
the summand of $R$ in Alexander grading $s$ and Maslov grading $d$,
$M_{s,d}=R_{s+1,d+1}\oplus R_{s,d}$, endowed with
the differential
$$D(x,y)=(\partial x, (U_b-U_c)x -\partial y).$$
This is quasi-isomorphic to the complex $\CFKm(\Knot_0)/U_b-U_c$.

\begin{theorem}
  \label{thm:StrongForm}
  Let $\Knot_+$ be a knot or link with a distinguished positive
  crossing, and let $U_a$ and $U_b$ be variables corresponding to
  the two out-going edges. There is a chain map $f\colon
  \CFLm(\Knot_+)\longrightarrow \CFLm(\Knot_-)$ whose mapping cone $E$
  is quasi-isomorphic to the mapping cone $M$ of the chain map
  $$U_{b}-U_c \colon \CFLm(\Knot_0)\longrightarrow \CFLm(\Knot_0).$$
  In the case where both strands at $\Knot_+$ belong to the same component
  of the knot, the quasi-isomorphism respects the bigrading, while
  in the case where the strands belong to different components,
  $E_{m,s}=M_{m,s-1}$.
\end{theorem}

We will give two proofs of the above theorem. But first, we show 
that it implies Theorem~\ref{thm:StrongForm}.

\begin{proof}[Theorem~\ref{thm:StrongForm} $\Rightarrow$ Theorem~\ref{thm:GradedSkein}]
  Suppose that $\Knot_+$ is connected. In this case, 
  skein exact sequence for $\HFKm$ follows immediately from Theorem~\ref{thm:StrongForm}, and the
  long exact sequence associated to a mapping cone.
  Consider next the case of $\CFKa$.
  Then, we have that
  \begin{align*}
  \CFKa(\Knot_+)&=\CFKm(\Knot_+)/U_a=0 \\
  \CFKa(\Knot_-)&=\CFKm(\Knot_-)/U_a=0 \\
  \CFKa(\Knot_0)&=\CFKm(\Knot_0)/U_a=U_b=0.
  \end{align*}
  Specializing our exact triangle to $U_a=0$, we obtain
  a long exact sequence connecting connect $\HFKa(\Knot_+)$,
  $\HFKa(\Knot_-)$, and $H_*(\CFKm(\Knot_0)/(U_a=0, U_a=U_b))$. Since the
  variables $U_a$ and $U_b$ correspond to basepoints $O_a$ and $O_b$
  correspond to two different components of $\Knot_0$, we can identify
  the latter homology group with $\HFKm(\Knot_0)$ as desired.
  
  Suppose that $\Knot_+$ consists of two components both of which
  project to the distinguished crossing. This time, we have
  \begin{align*}
    \CFKa(\Knot_+)&= \CFKm(\Knot_+)/(U_a=0,U_b=0) \\
    \CFKa(\Knot_-)&= \CFKm(\Knot_-)/(U_a=0, U_b=0) \\
    \CFKa(\Knot_0)&= \CFKm(\Knot_0)/U_a=0.
  \end{align*}
  Specializing our exact triangle to $U_a=0=U_b$, we connect
  $\HFKa(\Knot_+)$, $\HFKa(\Knot_-)$, and the homology of the mapping
  cone of $U_a-U_b$ on $\CFKm(\Knot_0)/(U_a=0=U_b)$.  In $\Knot_0$,
  since $a$ and $b$ belong to the same strand, multiplication by
  $U_a$ is homotopic to multiplication to $U_b$ (this follows from
  general properties of stabilization, cf.
  \cite[Section~\ref{Links:subsec:SimpleStabilizations}]{Links},
  \cite[Proposition~\ref{MOS:prop:ExtraBasepoints}]{MOS}, see
  also~\cite[Lemma~\ref{lemma:ChainHomotopiesZ}]{MOST} for a proof
  using grid diagrams). Thus, the mapping cone of
  $U_a-U_b$ on $\CFKm(\Knot_0)/(U_a=0=U_b)$ is quasi-isomorphic to
  the tensor product of $\CFKm(\Knot_0)$ with $V$.
  Similarly, for $\CFKm$, we have a
  triangle connecting $\CFKm(\Knot_+)$, $\CFKm(\Knot_-)$ and
  $\CFKm(\Knot_0)/(U_a-U_b)$.  Again, since 
  $O_a$ and $O_b$ correspond to the same component,
  $U_a-U_b$ is null-homotopic, so the latter complex is
  quasi-isomorphic to $\CFKm(\Knot_0)\otimes W$.

  Grading shifts are straightforward to verify (see the first proof of
  Theorem~\ref{thm:StrongForm} for more discussion on this).

  The more general case where $\Knot_+$ consists of more components
  follows from the same reasoning as above, but a little bit of extra notation.
\end{proof}

We give two proofs of Theorem~\ref{thm:StrongForm}. The first is a
pseudo-holomorphic curves proof, and the second combinatorial proof
uses grid diagrams. 

\subsection{Holomorphic curves proof}

\begin{figure}[ht]
\mbox{\vbox{\epsfbox{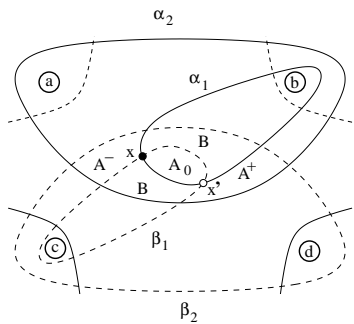}}}
\caption{\label{fig:ExactSequence}
{\bf{Exact triangle.}} Markings near a crossing used in
in the first proof of Theorem~\ref{thm:StrongForm}.}
\end{figure}

Our first proof of Theorem~\ref{thm:StrongForm} involves inspecting a
suitable Heegaard diagram, pictured in Figure~\ref{fig:ExactSequence}.
(This is also the route towards proving the exact triangle for singular
knots appearing in~\cite{Resolutions}.)

Draw a Heegaard diagram near a crossing as
shown in Figure~\ref{fig:ExactSequence}. In that picture, we have
distinguished circles $\alpha_1$ and $\beta_1$ which meet in two
points $x$ and $x'$.   This diagram is marked with $\Os=\{O_1,...,O_n\}$,
and $\Xs=\{A^-,A^0\}\cup \Xs_0$, where $\Xs_0=\{X_1,...,X_{n-2}\}$.
Alternatively, if we leave in $\alpha_1$ and $\beta_1$, and use
$\Xs=\Xs_0\cup \{A^-,A^0\}$, we obtain a Heegaard diagram for $\Km$.
Leaving in $\alpha_1$ and $\beta_1$, and using $\Xs=\Xs_0\cup
\{A^0\cup A^+\}$, we obtain a Heegaard diagram for the knot with
positive crossing $\Kp$.  Finally, using $\Xs$ as the union of $\Xs_0$
and the two regions marked by $B$, we obtain a Heegaard diagram for
the smoothing $\Knot_0$ of the crossing.

Note that there are four circles of type $\Os$ in the picture,
two of which correspond to the outgoing edges $a$ and $b$, and two
of which corresponding to the incoming ones $c$ and $d$.

Clearly, $\CFKm(\Km)$ has a subcomplex $X$ consisting of configurations
which contain the intersection point $x$, and a quotient complex $Y$. 
Thus, 
$\CFKm(\Km)$ can be thought of as the mapping cone of the map
$$D_{B}\colon Y \longrightarrow X$$
gotten by counting Maslov index one holomorphic disks representatives
of homology classes $\phi$ which contain exactly one of the regions
marked by $B$ (and hence neither of the regions marked $A^0$ or
$A^-$), and also none of the ones marked by the other $\Xs_0$,; i.e.
$$D_{B}(\x)=\sum_{\y\in{\mathfrak S}} \sum_{\left\{\phi\in\pi_2(\x,\y)
  \big| 
\begin{tiny}
\begin{array}{l}
\Mas(\phi)=1, \\
 X_i(\phi)=0~~~~~~~~~~\forall i=1,...,n-2 \\
B_1(\phi)+B_2(\phi)=1  
\end{array}
\end{tiny}
\right\}} \#\UnparModFlow(\phi) \cm
U_1^{O_1(\phi)}\cm...\cm U_n^{O_n(\phi)}\cm \y. $$
The understanding here is that the complex $X$ (and also $Y$) is
endowed with 
an induced differential
which counts holomorphic disks which do not cross any of
the four basepoints $A^0$, $A^-$, $B_1$, or $B_2$.
(There are, however, no constraints placed on the multiplicity
in $A^+$. It is not difficult to see, though, that the other constraints
imply that $A^+$ can be crossed at most once, and only for the differential within
$Y$).

Moreover,
$\CFKm(\Knot_0)$ has $Y$ as a subcomplex, with quotient $X$, and
hence, it can be thought of as the mapping cone of the map
$$D_{A^-}\colon X \longrightarrow Y,$$
defined by counting
flowlines which contain exactly one of the regions marked by $A^0$ or
$A^-$. 

Similarly, there is a subcomplex $X'$ of $\CFKm(\Knot_0)$
consisting of configurations which contain the intersection point
$x'$.  This has a quotient complex we denote by $Y'$. Moreover, $\Kp$
has a subcomplex isomorphic to $Y'$ and quotient complex isomorphic to
$X'$. Thus, we can think of $\CFKm(\Kp)$ as a mapping cone of
$$D'_{B}\colon X'\longrightarrow Y'$$
gotten by counting flowlines which contain exactly one of 
the regions marked by $B$, and neither of the regions
marked by $A^0$ or $A^+$.
Similarly, we can think of
$\CFKm(\Knot_0)$ as the mapping cone of
$$D_{A^+}\colon Y'\longrightarrow X',$$
which counts flowlines through exactly one of $A^0$ or $A^+$, and neither
of the regions marked by $B$.

There is an obvious isomorphism
$I\colon X \longrightarrow X'$, gotten by replacing the component $x$ by $x'$.
It is straightforward to verify that this is a chain map.

Consider the maps
\begin{eqnarray*}
D_{A^- B} \colon X\longrightarrow X &{\text{and}}&
D'_{A^+ B} \colon X'\longrightarrow X',
\end{eqnarray*}
where here $D_{A^- B}$ is defined by counting holomorphic disks
modulo translation in Maslov index one
homotopy classes $\phi$ with $X_j(\phi)\equiv 0$ $\forall j$, and
satisfying the addition conditions that
\begin{eqnarray*}
  A_0(\phi)+A^-(\phi)=1 &\text{and}& B_1(\phi)+B_2(\phi)=0;
\end{eqnarray*}
similarly, define $D_{A^+ B}$ to count holomorphic disks in homotopy
classes $\phi$ with
\begin{eqnarray*}
  A_0(\phi)+A^+(\phi)=1 &\text{and}& B_1(\phi)+B_2(\phi)=0;
\end{eqnarray*}

\begin{lemma}
  \label{lemma:HomotopyOne}
  The following relations hold:
  \begin{align*}
    D\circ D_{A^- B} + D_{A^- B}\circ D + D_B\circ D_A&=U_a+U_b-U_c-U_d \\
    D'\circ D_{A^+ B} + D_{A^+ B}\circ D' + D_A^+\circ D_B'&=U_a+U_b-U_c-U_d,
  \end{align*}
  where $D$ and $D'$ denote the differentials on $X$ and $X'$ respectively,
  and the right-hand-side represents multiplication by the scalar
  $(U_a+U_b-U_c-U_d)$ (thought of as an endomorphism of $X$ or $X'$).
  Informally, one can think of 
  $D_{A^- B}$ as furnishing a homotopy between $D_{B}\circ D_{A^-}$ and
  multiplication by $U_a + U_b - U_c - U_d$;
  and $D_{A^+ B}$ as furnishing a homotopy between $D_{A^+}\circ D'_{B}$
  and multiplication by $U_a+U_b-U_c-U_d$. 
\end{lemma}

\begin{proof}
        This is analogous to the proof (which uses Gromov's
        compactness theorem~\cite{Gromov}) that $\partial^2=0$ in Floer homology
        (cf.~\cite{HolDisk}, and~\cite{FOOO} for a general
        discussion).

        We consider the case of $A^- B$.  Look at ends of
        one-dimensional moduli spaces connecting $\x$ to $\y$ for
        homotopy classes $\phi$ which satisfy $A_0(\phi)+A^-(\phi)=1$
        and $B_1(\phi)+B_2(\phi)=1$.  These ends consist either of
        broken flowlines, or boundary degenerations. Broken flowlines
        are parameterized by pairs of homotopy classes of Maslov index one
        homotopy classes
        $\phi_1\in\ModFlow(\x,\x')$, $\phi_2\in(\x',\y)$ for some $\x'\in\Gen$.
        These can be partitioned into four cases: 
        \begin{itemize}
        \item $A_0(\phi_1)+A^+(\phi_1)=1$
        and $B_1(\phi_1)+B_2(\phi_1)=1$
        \item $A_0(\phi_1)+A^+(\phi_1)=0$
        and $B_1(\phi_1)+B_2(\phi_1)=1$
        \item $A_0(\phi_1)+A^+(\phi_1)=1$
        and $B_1(\phi_1)+B_2(\phi_1)=0$
        \item $A_0(\phi_1)+A^+(\phi_1)=0$
        and $B_1(\phi_1)+B_2(\phi_1)=0$.
        \end{itemize}
        The first types are counted in $D\circ D_{A^- B}$,
        the second by $D_{A^-}\circ D_{B}$, the third
        $D_{B}\circ D_{A^-}$, and the fourth in $D_{A^- B}\circ D$.

        The contributing boundary degenerations in the ends of this
        moduli space consist of Maslov index two holomorphic disks
        with boundary in $\Ta$ or $\Tb$, and which contain both one
        point from $\{A_0, A^-\}$ and one in $\{B_1,B_2\}$. There are
        four of these, one of which contains each of $O_a$, $O_b$,
        $O_c$, or $O_d$ respectively
        (compare~\cite[Lemma~\ref{Resolutions:lemma:CalcComposite}]{Resolutions}).
\end{proof}

We form now the chain complex $C$, given by the diagram:
\[
\begin{diagram}
X &\rTo^{D_{A^-}} & Y \\
\dTo^{D'_{B}\circ I} & \rdTo^{H} & \dTo_{I\circ D_{B}} \\
Y' & \rTo_{D_{A^+}} & X'
\end{diagram}
\]
where $H=-I\circ D_{A^-B}+ D'_{A^+ B}\circ I$.
In fact, the above diagram can be used to form a chain 
complex thanks to Lemma~\ref{lemma:HomotopyOne}. We denote this complex by $E$.

Clearly, the above complex has a subcomplex, corresponding to the
rightmost column, which is the mapping cone of
$-I\circ D_B\colon Y\longrightarrow X'$, which in turn is identified with
$\CFKm(\Knot_-)$, while its quotient complex is the mapping cone
of $D'_{B}\circ I \colon X \longrightarrow X$, which in
turn is identified with $\CFKm(\Knot_+)$.

Moreover, the bottom row is a subcomplex in turn is identified with
$\CFKm(\Knot_0)$; its quotient complex is the top row which also
is identified with $\CFKm(\Knot_0)$.

\begin{lemma}
  Under the identification of both rows of $E$ with $\CFKm(\Knot_0)$,
  the vertical map $D'_{B}\circ I + I\circ D_{B}$ is homotopic via
  $H$
  to multiplication by $U_b-U_c$.
\end{lemma}

\begin{proof}
  This follows along the lines of Lemma~\ref{lemma:HomotopyOne}.
  We can think of $I$ as the map $D_{A_0}$ gotten by counting holomorphic
  disks which cross $A_0$. Now the sum of vertical maps is induced
  by $D'_{B}\circ D_{A_0}+D_{A_0}\circ D_{B}$. Moreover, the map
  $D_{B A_0}$ induces a homotopy of this map with the count of all
  boundary degenerations containing both $A_0$ and $B$. This latter
  map is readily seen to correspond to multiplication by $U_b-U_c$.
\end{proof}

\begin{proof}[of Theorem~\ref{thm:StrongForm}]
  As we have seen, the complex $E$ is simultaneously identified with
  the mapping cone of a map $f\colon \CFLm(\Knot_+)\longrightarrow
  \CFKm(\Knot_-)$, and the mapping cone of a map
  $\CFKm(\Knot_0)\longrightarrow \CFKm(\Knot_0)$ which is chain
  homotopic to multiplication by $U_a-U_c$, which in turn is quasi-isomorphic
  to $\CFKm(\Knot_0)/U_b-U_c$, as desired.
  
  We turn our attention to gradings. Configurations $X$ and $Y$ inherit Maslov
  and Alexander gradings from either $\Knot_-$ or $\Knot_0$;
  similarly, configurations in $X'$ or $Y'$ inherit Maslov and Alexander
  gradings from either $\Knot_0$ or $\Knot_-$.  We
  assert that all the induced Maslov gradings coincide. This
  is clear since the Maslov grading of a given generator
  is independent of the placement of
  points of type $\Xs$, depending only on the placement
  of the $\Os$ (which coincide for all three links).
  
  Consider next the $N=M-2A$-gradings. After setting all $U_i=1$, we
  can isotope across the $\Os$ to obtain the same diagram for
  $\Knot_-$ and $\Knot_+$.  Thus the $N$-gradings of the two diagrams
  agree.
  Thus, it follows that (absolute)
  $A$-gradings for $\Knot_-$ and $\Knot_+$ coincide for all
  generators.

\begin{figure}[ht]
\mbox{\vbox{\epsfbox{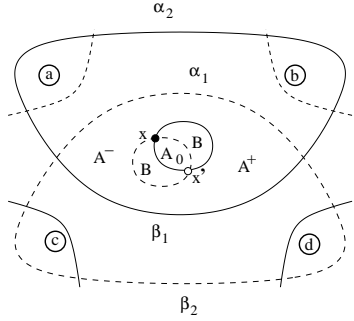}}}
\caption{\label{fig:AfterIsotopy}
{\bf{Isotoping across $O_b$ and $O_c$.}} 
In the complex where $U_b=U_c=1$, we can isotope across $O_b$ and $O_c$
without changing homology. The resulting picture is shown here.}
\end{figure}

Now consider the horizontal connecting homomorphism $f\colon
\CFLm(\Knot_+)\longrightarrow \CFKm(\Knot_-)$.  This map clearly
preserves both Alexander and Maslov gradings.  For example, we can
view the restriction of $f$ to the subcomplex $Y'$. The component of
the connecting homomorphism is gotten by counting holomorphic disks
which cross exactly one of $A_0$ or $A^+$, a map which simultaneously
drops $A$ and $M$-degree by one, post-composed by the inverse of $I$,
which simultaneously raises both of these degrees by one.
  
  Using the bigrading grading on $E$ for which the projection map
  $\pi$ respects bigradings so that $\pi\colon E_{m,s} \longrightarrow
  \CFKm_{m-1}(\Knot_+,s)$, we then have a bigraded identification
  $$E_{m,s}\cong X'_{m,s} \oplus Y'_{m,s}\oplus Y_{m+1,s}\oplus
  X_{m+1,s},$$ where all identifications are made using the diagram for
  $\Knot_+$.  In fact, under the bigraded isomorphism with the mapping
  cone of $U_b-U_c$, we have that $$E_{m,s}\cong X_{m+1,s+1}\oplus
  Y'_{m,s}\oplus Y_{m+1,s}\oplus X'_{m+1,s+1},$$ where once
  again all bigradings are induced from $\Knot_+$. 
  
  Suppose now that the two strands in $\Knot_+$ belong to the same
  component.  In this case, we claim $\Knot_0$ and $\Knot_+$ induce
  the same the bigrading on $X$. To see this, observe that in the
  $U_i\equiv 1$ complex, a generator $\x$ for $\CFKm(\Knot_+)$, when
  thought of as an element of $H_*(T^{n-1})$ has grading one greater than the
  same same element thought of as a generator for for $\CFKm(\Knot_0)/\{U_i\equiv 1\}$.
  Thus, we have a bigraded identification
  $$E_{m,s}\cong R_{m+1,s+1}\oplus R_{m,s}.$$
  Similarly, if the two
  strands in $\Knot_+$ belong to different components, then the
  Alexander grading of a generator from $X$ thought of as represented
  in $\Knot_0$ is one less than its Alexander grading thought of as a
  represented in $\Knot_+$; hence we have that
  $$E_{m,s}\cong R_{m+1,s}\oplus R_{m,s-1}.$$
\end{proof}

\subsection{Proof using grid diagrams}

\begin{figure}[ht]
\mbox{\vbox{\epsfbox{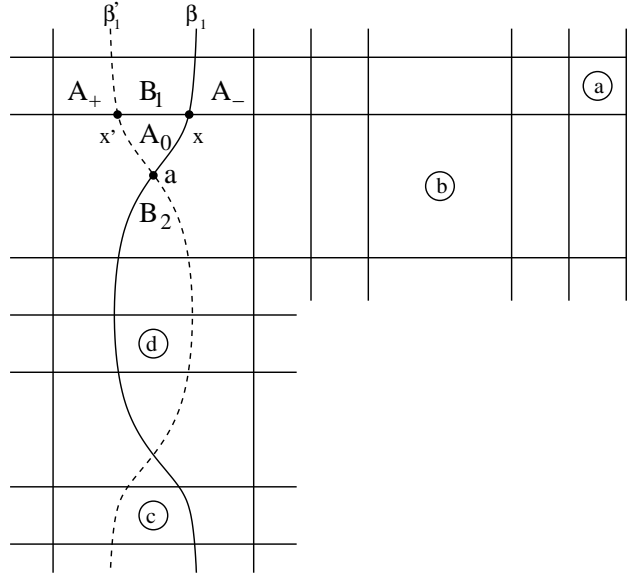}}}
\caption{\label{fig:GridProof}
{\bf{Grid diagram for the skein relation.}}  Grid diagrams for a proof
of Theorem~\ref{thm:StrongForm}.}
\end{figure}

A disadvantage of the grid diagrams proof is that there is no grid
diagram which represents all three knots $\Knot_+$, $\Knot_0$, and
$\Knot_0$; instead, one has to move the grid diagram, as pictured in
Figure~\ref{fig:GridProof}. By way of explanation, we have the
Heeegaard torus, equipped with $\Os$, $\Xs$, and horizontal circles
$\alphas$.  We have two possible sets of vertical circles $\betas$ and
$\betas'$, which differ only in the choice of the first circle
(i.e. $\beta_i=\beta_i'$ for $i>1$), $\beta_1$ and $\beta_1'$. We call
the two grid diagrams $G$ and $G'$. The circles $\beta_1$ and $\beta'$
meet in two points, one of which is labelled $a$ as shown in the
picture. Note that there is a small triangle bounded by an arc in an
$\alpha$-circle, an arc in $\beta_1$, and an arc in $\beta_1'$, which
contains $A_0$ in its interior, and whose three vertices are $x$ (on
$\beta_1$), $x'$ (on $\beta_1'$) and $a$.

In the present context, the complex $X$ is generated by generators
$\Gen(G)$ which contain the intersection point $x$, while $X'$ is
generated by generators $\Gen(G')$ which contain the intersection
point $x'$. These are made into complexes by counting rectangles which
are disjoint from $B_1$ and $B_2$.  There are also complexes
$Y$ and $Y'$ defined using the complementary sets of generators.
As before, we have maps
\begin{align*}
D_{A^-}&\colon X \longrightarrow Y \\
D_B &\colon Y \longrightarrow X \\
D_B'&\colon X' \longrightarrow Y'\\
D'_{A^+}&\colon Y'\longrightarrow X'
\end{align*}
The first two of these maps is defined by by counting rectangles in
$G$: the first counts rectangles which contain one of $A_0$ or $A_-$,
the second counts rectangles which contain one of $B_1$ or $B_2$.  The
second two use rectangles in the diagram $G'$, the first counts
rectangles containing $B_1$ or $B_2$, and the second counts rectangles
containing $A_0$ or $A_+$.
There is also an identification $I\colon X \longrightarrow X'$
gotten by moving the intersection point $x$ to $x'$.

\begin{lemma}
        \label{lemma:GridComplex}
        We have the following relations:
        \begin{align*}
        D_B\circ D_{A_-}|_{X}&=U_a+U_b-U_c-U_d \\
        D_{A^+}\circ D'_{B}|_{X'}&=U_a+U_b-U_c-U_d
        \end{align*}
\end{lemma}

\begin{proof}
        This is essentially a repetition of Lemma~\ref{lemma:HomotopyOne},
        except, of course, that for grid diagrams, arguments like Gromov's
        compactness theorem can be formulated in purely combinatorial  terms
        (cf.~\cite[Proposition~\ref{MOST:prop:DSquaredZero}]{MOST}).

        We start with the first equation. Observe that the sum of maps
        $D_B\circ D_{A^-}+D_{A^-}\circ D_B$ counts polgons obtained by
        juxtaposing two rectangles, one of which contains $A_0$ or
        $A_-$, and the other contains $B_1$ or $B_2$.  These polygons
        cancel in pairs, except for those annuli of length or width
        equal to one, which contain both $A_0$ or $A_-$ and $B_1$ or
        $B_2$, of which there are four, contributing
        $U_a+U_b-U_c-U_d$. In principle, there might be alternative
        decompositions consisting of a rectangle without either $A_0$,
        $A_-$, $B_1$, or $B_2$ (and the other must have one of each
        pair).  But there are no such alternative decompositions with
        initial point at $X$. In a similar vein, it is straightforward
        to see that $D_B$ annihilates $X$. 

        The second relation follows similarly.
\end{proof}

According to Lemma~\ref{lemma:GridComplex}, the following diagram
represents a chain complex:
\begin{equation}
        \label{eq:GridComplex}
\begin{CD}
X @>{D_{A^-}}>>  Y \\
@V{D'_{B}\circ I}VV  @VV{I\circ D_{B}}V \\
Y' @>{D_{A^+}}>> X'
\end{CD}
\end{equation}

Let $C$ be the chain complex for the resolution $\Knot_0$ using
$\beta_1$ in Figure~\ref{fig:GridProof}, and $C'$ be the chain complex
for $\Knot_0$ using $\beta_1'$.  An explicit chain homtopy equivalence
$$\Phi\colon C'\longrightarrow C$$ is constructed
in~\cite[Subsection~\ref{MOST:subsec:Commutation}]{MOST}.  This map is
defined by counting pentagons.  More precisely, given $\x\in \Gen(G)$
and $\y\in\Gen(H)$, we let $\Pent(\x,\y)$ denote the space of embedded
pentagons with the following properties. This space is empty unless
$\x$ and $\y$ coincide at $n-2$ points. An element of
$\Pent_{\beta'\beta}(\x,\y)$ is an embedded disk in $\Torus$, whose
boundary consists of five arcs, each contained in horizontal or
vertical circles. Moreover, under the orientation induced on the
boundary of $p$, we start at the $\beta_1'$-component of $\x$,
traverse the arc of a horizontal circle, meet its corresponding
component of $\y$, proceed to an arc of a vertical circle, meet the
corresponding component of $\x$, continue through another horizontal
circle, meet the component of $\y$ contained in $\beta_1$, proceed to
an arc in $\beta_1$ until we meet the intersection point $a\in
\beta_1\cap\beta'$, and finally, traverse an arc in $\beta_1'$ until
we arrive back at the initial component of $\x$.  Finally, all the
angles here are required to be less than straight angles. The space of
empty pentagons $p\in \Pent_{\beta'\beta}(\x,\y)$ with $\x \cap
\Interior(p) = \emptyset$, is denoted $\EmptyPent_{\beta'\beta}$.

Given $\x\in\Gen(G')$, define
\begin{align*}
\Phi(\x) =& \sum_{\y\in \Gen(G')}\,\,
    \sum_{p\in\EmptyPent_{\beta'\beta}(\x,\y)}\!\!  U_1^{O_1(p)}\cdots
    U_n^{O_n(p)} \cm \y
\in \Cm(G).
\end{align*}
It is elementary to see that the above map induces a chain homotopy
equivalence~\cite[Proposition~\ref{MOST:prop:Commute}]{MOST}.

\begin{lemma}
  \label{lemma:VerticalMap}
        Under the natural identification of the two rows in the
        above complex with chain complexes for $\CFKm(\Knot_-)$,
        the two vertical maps add up to multiplication by $U_b-U_c$.
\end{lemma}

\begin{proof}
        Let $C$ be the chain complex for $\CFKm(K_0)$
        appearing in the top row of Equation~\eqref{eq:GridComplex},
        and let
        $C'$ be the complex for $\CFKm(K_0)$ appearing in the bottom
        row. Thus, $C$ belongs to the grid diagram $G$, while $C'$ belongs
        to the grid diagram $G'$.
        The sum of the two vertical maps can be viewed as a chain map
        $V\colon C \longrightarrow C'$,
        where 
        $$V=D'_{B}\circ I \circ \Pi_X + I \circ D_B \circ \Pi_Y.$$
        
        We have a map $\Phi\colon C'\longrightarrow C$ defined by
        counting pentagons, as above.
        
        We will find it convenient to extend the maps $D_B$ and $D'_B$
        earlier to maps
        \begin{eqnarray*}
        \eD_B\colon C \longrightarrow C &{\text{and}}&\eD'_B\colon
        C'\longrightarrow C',
        \end{eqnarray*} gotten by counting rectangles which contain
        both $B_1$ and $B_2$. Similarly, we have a map $\eD_{A_0}$ and
        $\eD_{A_0}'$ defined by counting rectangles which contain
        $A_0$.

        We claim that
        \begin{equation}
        \label{eq:LeftColumn}
        \Phi\circ D'_{B}\circ I \circ
        \Pi_X={\eD}_{B}\circ \eD_{A_0}.
        \end{equation}
        This is seen as follows. The
        composite $\Phi\circ {D_B}\circ I \circ \Pi_X$ is a count of
        polygons, which are gotten by juxtaposing an empty rectangle
        starting at an intersection point from $X$, followed by an
        empty pentagon. By filling in the small triangle containing
        $A_0$, we obtain a one-to-one correspondence between these
        polgons, and polygons obtained in the following way:
        \begin{enumerate}
        \item 
        juxtapositions of two rectangles, the first of which contains $A_0$
        and the second of which contains $B_2$
        \item 
        juxtapositions of rectangles, the first of which is empty,
        and the second of which contains both $A_0$ and $B_1$
        \item 
        the column {\em in $G'$} through both $A_0$ and and $B_1$.
      \end{enumerate}
      The first term contibutes ${\widetilde D}_{B_2}\circ {\widetilde
        D}_{A_0}$.  Decomposing the polygon in an alternative way, we
      see that the sum of the second two terms contributes
      $\eD_{B_1}\circ \eD_{A_0}$. (Note that the column in $G'$ through 
      $A_0$ and $B_1$ contributes $U_c$, which is the same as the contribution
      the column in $G$ of the column through $A_0$ and $B_1$.)
        
        Similarly, we claim that        
        \begin{equation}
        \label{eq:RightColumn}
        \Phi \circ I \circ D_{B} \circ \Pi_Y=\eD_{A_0}\circ \eD_B.
        \end{equation}
        This follows more directly than Equation~\eqref{eq:LeftColumn}.
        Filling in the small triangle containing
        $A_0$, we obtain a one-to-one correspondence between the polgons
        counted on both sides.
        
        Adding Equations~\eqref{eq:LeftColumn}
        and~\eqref{eq:RightColumn}, we conclude that $\Phi\circ V=\eD_{A_0}\circ
        \eD_B + \eD_B\circ \eD_{A_0}$. The same arguments from Lemma~\ref{lemma:GridComplex}
        show that the following relation holds:
        $$D \circ \eD_{A_0 B}+\eD_{A_0 B} \circ D = \eD_{B}\circ \eD_{A_0} + \eD_{A_0}\circ D_{B} + U_b-U_c;$$
        i.e. $\Phi\circ V$ is homotopic to multiplication by $U_b-U_c$, as desired.

\begin{figure}[ht]
\mbox{\vbox{\epsfbox{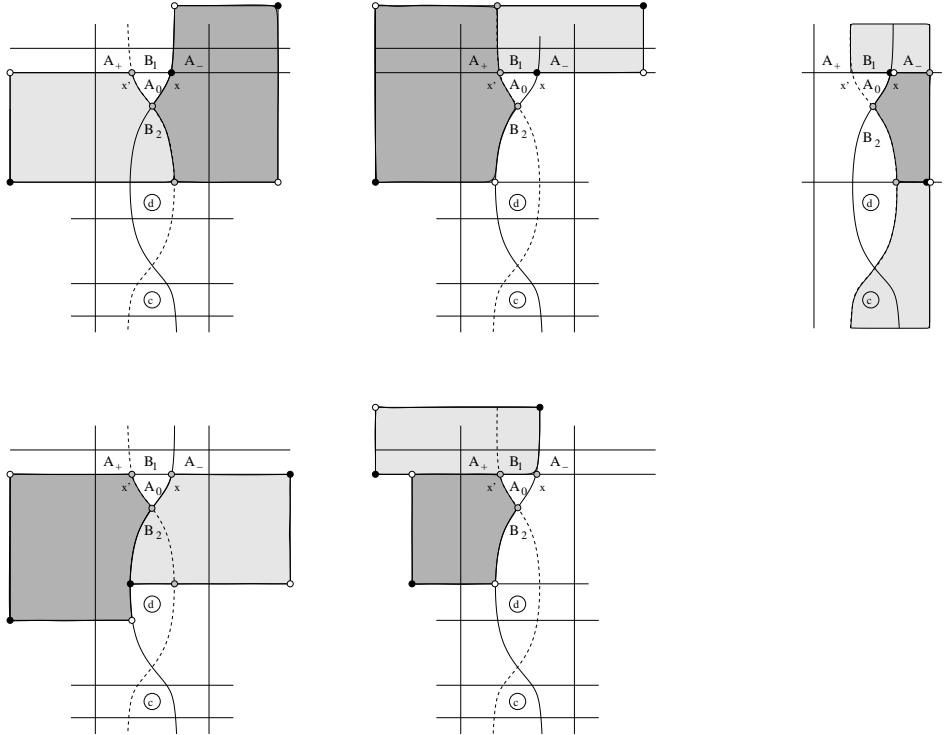}}}
\caption{\label{fig:Cases}
{\bf{Proof of Lemma~\ref{lemma:VerticalMap}.}} 
        In the top row, we have illustrated the polygons
        contributing to $\Phi\circ D_{B}\circ I \circ
        \Pi_X$, in the second, we have illustrated the ones 
        contributing to $\Phi \circ I \circ D_{B} \circ \Pi_Y$.
        The lighter polgon is always composed first; the dark dots
        represent initial points. Filling in the triangle labelled
        by $A_0$ gives a one-to-one correspondence between these polygons
        and those counted in $\eD_B \circ \eD_{A_0}+\eD_{A_0}\circ \eD_B$.}
\end{figure}
\end{proof}

\begin{proof}[Grid diagram proof of Theorem~\ref{thm:StrongForm}.]
  The proof follows from inspecting the complex from
  Equation~\eqref{eq:GridComplex}, which we refer to now as $E$.
  Again, the two vertical columns correspond to $\CFKm(\Knot_-)$ and
  $\CFKm(\Knot_+)$ respectively, and so the horizontal maps add
  up to give the stated map $f$. Moreover, according to 
  Lemma~\ref{lemma:VerticalMap}, $E$ is
  quasi-isomorphic to the mapping cone of 
  $$U_b-U_c\colon \CFKm(\Knot_0)\longrightarrow \CFKm(\Knot_0).$$
  Gradings can be traced through exactly as in the earlier proof.
\end{proof}

\bibliographystyle{plain}
\bibliography{biblio}

\end{document}